\documentclass[12pt]{article}

\newcommand{\RED}[1]{{\color{red}{#1}}}
\newcommand{\BLUE}[1]{{\color{blue}{#1}}}

\usepackage{amssymb}
\usepackage{amsfonts}
\usepackage{latexsym}
\usepackage{graphicx}

\usepackage[usenames]{color}
\usepackage{amsmath,amsthm}

\usepackage{color}

\newcommand{\fA}{\mathfrak{A}}

\newcommand{\fG}{\mathfrak{G}}

\newcommand{\fN}{\mathfrak{N}}
\newcommand{\fX}{\mathfrak{X}}

\newcommand{\fF}{\mathfrak{F}}

\newcommand{\fS}{\mathfrak{S}}

\newcommand{\fZ}{\mathfrak{Z}}

\newcommand{\qn}{\mathrm{q}{\rm Max}\mbox{-}\lhd}

\newcommand{\mda}{{\rm Max}\mbox{-}\lhd\fA}

\newcommand{\QQQ}{\mbox{\footnotesize \rm Q}}
\newcommand{\EEE}{\mbox{\footnotesize \rm E}}
\newcommand{\PPP}{\mbox{\footnotesize \rm P}}
\newcommand{\HHH}{\mbox{\footnotesize \rm H}}
\newcommand{\LLL}{\mbox{\footnotesize \rm L}}

\newcommand{\ignore}[1]{}

\newcommand{\sn}{\mathrm{sn}}

\newcommand{\beqn}{\begin{eqnarray*}}
\newcommand{\eeqn}{\end{eqnarray*}}

\newcommand{\aut}{\mathrm{aut}}

\newcommand{\Z}{\mathbb{Z}}

\newcommand{\R}{\mathbb{R}}

\newcommand{\C}{\mathbb{C}}

\newcommand{\N}{\mathbb{N}}

\newtheorem{theorem}{Theorem}[section]
\newtheorem{lemma}[theorem]{Lemma}
\newtheorem{definition}[theorem]{Definition}
\newtheorem{corollary}[theorem]{Corollary}

\newtheorem{remark}[theorem]{Remark}

\renewcommand{\epsilon}{\varepsilon}
\newcommand{\ONE}{{\bf 1}}

{\begingroup

  \begin{enumerate}}
  {\end{enumerate}\endgroup}

\title{Chain Conditions for Abelian Normal Subgroups}

\author{Falih A.M. Aldosray\\
        Department of Mathematics\\
	Umm Al-Qura University\\
	PO Box 56199 \\ 
	Makkah, Saudi Arabia 
               \and
        Ian Stewart \\
        Mathematics Institute \\
        University of Warwick \\
        Coventry CV4 7AL \\
        United Kingdom
}

\begin{document}

\maketitle

\begin{abstract}

We investigate various chain conditions for abelian normal subgroups of
groups, following Robinson \cite{R68}. We focus on the
maximal condition for abelian normal subgroups, an exceptional
case excluded from his results except for some counterexamples.
In particular we show that a hypercentral group satisfies the
maximal condition for abelian normal subgroups if and only if it is 
is finitely generated nilpotent.
\end{abstract}

\section{Introduction}

Robinson \cite{R68} studied chain conditions Max-$\lhd^\alpha \fX$
and $\lhd^\alpha \fX$-Max for groups. Here $\fX$ is a class of groups in the
sense of Philip Hall and $\alpha$ is an ordinal. We write $H \lhd^\alpha G$
when $H$ is an $\alpha$-step ascendant subgroup of $G$. If $\alpha = n$ is finite,
$H$  is in $n$-step subnormal subgroup of $G$.

A group $G$ lies in Max-$\lhd^\alpha \fX$ if the set of all $\lhd^\alpha$
subgroups of $G$ that lie in $\fX$ satisfies the maximal condition:
every ascending chain terminates after finitely many steps. 
A group $G$ lies in $\lhd^\alpha \fX$-Max if every $\lhd^\alpha$
subgroup of $G$ that lies in $\fX$ satisfies the maximal condition
for subgroups --- equivalently, every subgroup is finitely generated. 
There are also analogous classes
Min-$\lhd^\alpha \fX$ and $\lhd^\alpha \fX$-Min, where the maximal
condition is replaced by the minimal condition.

The adjective `Noetherian', suitably embellished, is also used to
indicate a maximal condition. Similarly, `Artinian' is often used to
indicate a minimal condition. These terms originated in ring theory.
We do not consider Artinian conditions in this paper.

The notation Max denotes the class of Noetherian groups: those that satisfy the
maximal condition for subgroups. Similarly, Min denotes the class of Artinian groups: those that satisfy the
minimal condition for subgroups.

Robinson \cite{R68} focuses mainly on the cases when $\fX$
is one of:

\qquad $\fA$ : abelian groups.

\qquad $\fN$ : nilpotent groups.

\qquad $\fS$ : soluble (solvable) groups. This class is also denoted
by $\EEE\fA$.

\noindent
He proved (Corollary 2) that if $\alpha \geq 2$ then the five classes
\[
{\rm Max}-\lhd^\alpha \fN \quad {\rm Max}-\lhd^\alpha \fS \quad
\lhd^\alpha \fA-{\rm Max} \quad \lhd^\alpha \fN-{\rm Max}\quad\lhd^\alpha \fS-{\rm Max}
\]
coincide. Moreover, the soluble radical $\sigma(G)$ of such a group (the join of all
soluble normal subgroups)
is always polycyclic, and contains every $\lhd^{\alpha-1} \fS$ subgroup of $G$.
In particular, $\sigma(G)$ is finitely generated (and satisfies Max).
He also investigated analogous issues for Min.

Section 4 of his paper discusses various counterexamples to these results
in the exceptional case $\alpha = 1$. We make further observations about
this exceptional case. In particular, we prove that
every hypercentral group that lies in ${\rm Max}-\lhd\fA$ is
finitely generated nilpotent.

\subsection{Notation and Terminology}
\label{S:NT}

Our notation and terminology for groups is standard, except that
we use a minor variant of the notation for
group classes and closure operations, which were
introduced by Hall~\cite{pH59,pH63}, and are described in detail in
Robinson~\cite[Section 1.1]{R70a}.
The main differences are that we write closure operations
as small Roman capitals, such as $\EEE$, instead of bold italic \textbf{\textit{E}}.
Also we use $\EEE$, not $\PPP$ or \textbf{\textit{P}}, for closure under extensions, and
$\QQQ$, not $\HHH$ or \textbf{\textit{H}}, for closure under quotients.
Additional notation is stated when required.
For clarity, we
briefly recall some of the basics.

If $x, g \in G$ then the {\em conjugate} of $x$ by $g$ is written
\[
x^g = g^{-1}xg
\]
The {\em commutator} of $x$ and $g$ is
\[
[x,g] = x^{-1}g^{-1}xg = x^{-1}x^g
\]
We denote the identity element of a group by $1$ unless otherwise stated.
The one-element group is denoted by ${\bf 1}$.
We write $H \leq G$ if $H$ is a subgroup of $G$, and
$H \lhd G$ if $H$ is a normal subgroup of $G$.

The subgroup generated by $X \subseteq G$ is denoted by $\langle X \rangle$.
If $X, Y \subseteq G$ then their {\em commutator subgroup} is
\[
[X,Y] = \langle [x,y] : x\in X, y \in Y \rangle
\] 
If $X = \{x\}$ we abbreviate this to $[x,Y]$, and similarly 
if $Y = \{y\}$. We do not use the convention if both $X$ and $Y$ 
are singletons since this would make $[x,y]$ ambiguous.

The {\em derived group} of $G$ is
\[
G' = \langle [x,y] : x,y, \in G \rangle
\]
and the {\em derived series} comprises the subgroups $G^{(m)}$
defined inductively by 
\[
G^{(0)} = G \qquad G^{(m+1)}=(G^{(m)})'
\]
so $G'=G^{(1)}$.
The group $G$ is {\em soluble} of {\em derived length} $m$ 
if $G^{(m)}=\ONE$ for some finite $m$.

The {\em centraliser} of a subgroup $H \leq G$ is
\[
C_G(H) = \{ x \in G : xh = hx\ \ \forall h \in H \}
\]
The {\em centre} of $G$ is $C_G(G)$, denoted $\zeta_1(G)$, and the
{\em upper central series} comprises the subgroups $\zeta_m(G)$
defined inductively by
\[
\zeta_{m+1}(G)/\zeta_m(G) = \zeta_1(G/\zeta_m(G))
\]
The {\em lower central series} comprises the subgroups $\gamma_m(G)$,
where $\gamma_1(G) = G'$ and
\[
\gamma_{m+1}(G) = [G,\gamma_m(G)]
\]
The group $G$ is {\em nilpotent} of {\em class} $m$ if $\zeta_m(G)=G$ for some finite $m$.
Equivalently, $\gamma_m(G)=\ONE$.

A group is {\em Noetherian} if it satisfies the ascending chain condition on normal subgroups.
The class of Noetherian groups is denoted by Max-$\lhd$. Following
\cite{AS26}, by analogy with the Lie algebra setting \cite{AS19},
a group $G$ is {\em quasi-Noetherian} if, given any ascending chain 
$N_1 \subseteq N_2 \subseteq \cdots$ of normal subgroups of $G$,
then there exists $r > 0$ such that $[G^{(r)},N_\infty] \subseteq N_r$,
where $N_\infty = \bigcup_i N_i$. The class of quasi-Noetherian groups 
is denoted by $\qn$. It is easy to prove that every Noetherian or soluble group 
is quasi-Noetherian, but there are others (see \cite{AS26}).

The {\em Fitting radical} (Fitting subgroup, nilpotent radical) $\nu(G)$ 
is the product of all nilpotent normal 
subgroups of $G$. It is locally nilpotent.

The {\em Hirsch--Plotkin radical} $\rho(G)$ is the product of all locally 
nilpotent normal subgroups of $G$. Hirsch \cite{H55} and Plotkin \cite{P54}
proved that it is the unique maximal
locally nilpotent normal subgroup; see also \cite[Section 2.3]{R70a}.

The {\em soluble radical} $\sigma(G)$ is the product of all soluble normal subgroups of $G$.
It is locally soluble.

\section{Maximal Condition for Abelian Normal Subgroups}

By analogy with~\cite[Chapter 9]{AS74} and~\cite[Section 4]{AS19},
and following Robinson \cite{R68},
we define the class $\mda$ to consist of those
groups satisfying the maximal condition for abelian normal subgroups.
We establish some properties related to hypercentrality for use in Section \ref{S:HGMFGN}.

\begin{lemma}
\label{L:MC1.1}
Let $R$ be the join of a collection of $\fF\cap\EEE\fA$ normal subgroups
of $G$ and suppose that $R$ is infinite. Then $R$ contains
an infinite hypercentral normal subgroup $H$ of $G$ with $H = \zeta_\omega(H)$.
\end{lemma}
\begin{proof}
Evidently we can find a strictly ascending chain
\[
\ONE = R_0 \subsetneq R_1 \subsetneq R_2 \subsetneq \cdots
\]
of  $\fF\cap\EEE\fA$ normal subgroups of $G$ such that $R_i \subsetneq R$
for all $i$. We can refine this chain to one with abelian factors, say
\[
\ONE = N_0 \subsetneq N_1 \subsetneq N_2 \subsetneq \cdots
\]
where $N_i \lhd G$ and $N_{i+1}/N_i \in \fA$ for all $i$. Let $N = \bigcup_{j=0}^\infty N_j$.
Then $N \lhd G$ and $N$ is infinite. 

Let $H_1 =N_1 \lhd G$, so $H_1 \in \fA$. Inductively, suppose we have defined
$H_i$ such that $H_i \in \fF\cap\EEE\fN$, with $H_i \lhd G$ and $H_i \subseteq N$. 
Since $H_i \in \fF$ we have $N/C_N(H_i) \in \fF$, so that $C_N(H_i) \not\in \fF$,
whence $C_N(H_i) \neq H_i$. Thus we can find $n_i$ minimal with respect to
$C_{n_i} = C_N(H_i) \cap N_{n_i} \not\subseteq H_i$, because
$C_N(H_i) = \bigcup N_j \cap C_N(H_i)$ and $H_i \in \fF$. Clearly
$n_i \neq 0$ and $C_{n_i}' \subseteq C_N(H_i) \cap N_{n_i-1}$,
because $N_{n_i}' \subseteq N_{n_i-1}$ and $n_i$ is minimal. Therefore
$\gamma_3(C_{n_i}) = \ONE$.

Define $H_{i+1} = C_{n_i}H_i$. Then $H_{i+1} \in \fF\cap\fN$ and
$H_{i+1} \subseteq N$. Now
$H_i \lhd G$ implies that $C_G(H_i) \lhd G$, and $N \lhd G$. Therefore $C_N(H_i) = C_G(H_i) \cap N \lhd G$, so
$C_{n_i} = C_N(H_i) \cap N_{n_i} \lhd G$. Therefore $H_{i+1} \lhd G$
and $H_i \subsetneq H_{i+1}$. Proceeding inductively, define
$H_i$ for all $i$ so that  $H_i \lhd G$, $H_i \in \fF\cap\fN$, and $H_i \subsetneq H_{i+1}
\subseteq N$. Evidently $H_{i+1}' \subseteq H_i$ and $H_{i+1}$ is centralised by
$C_{n_i+k}$ for all $k \geq 1$. Now let $H = \bigcup_{i=1}^\infty H_i$,
so that $H$ is infinite and $H \lhd G$. Clearly $H \langle H_{i+1}, 
C_{n_i+1}, C_{n_i+2}, \ldots\rangle$ for all $i = 1,2,\ldots$.Therefore
$[H,H_{i+1}] \subseteq H_{i+1}' \subseteq H_i$,
whence $H_i \subseteq \zeta_i(H)$ for all $i$. Thus $H = \zeta_\omega (H)$.
Finally, $H \subseteq N \subseteq R$ and the lemma is proved.
\end{proof}

Recall that a subgroup $H$ of a group $G$ is {\em self-centralising}
if $C_G(H) = H$. This implies in particular that $H$ is abelian.
Recall that a subgroup $H$ of $G$ is {\em characteristic} if it is
invariant under all automorphisms of $G$. We write $H\  \mathrm{ch}\ G$.
Characteristic subgroups are normal (consider inner automorphisms).
Moreover, $H\  \mathrm{ch}\ K \lhd G$ implies $H  \lhd G$.

\begin{lemma}
\label{L:MC1.2}

\noindent{\rm (a)}
A maximal abelian normal subgroup of a hypercentral group
is self-centralising.

\noindent{\rm (b)}
If $G$ is any group, then $\sigma(G) \in \acute{\EEE}\fA$, where
$\acute{\EEE}\fA$ is the class of all groups with an ascending abelian series.

\noindent{\rm (c)}
Let $H \lhd G$ and let $\Delta$ be any of the relations $\leq, \lhd, \lhd^k, \sn, \mathrm{asc}$.
If $G \in  \acute{\EEE}\fA$ and $H$ contains every $\fN_2$ $\Delta$-subgroup of $G$
then $C_G(H) \subseteq H$.

\noindent{\rm (d)}
Suppose that $H \lhd G$ and $H$ contains every $\fN_2$ normal subgroup of $G$.
If
\[
K = \langle N : N \lhd G\ \&\ N \in \EEE\fA \cup \fZ \rangle
\]
then
\[
K \cap C_G(H) = \zeta_1(H) \subseteq \zeta_1(H\cap K)
\]
In particular, for any group $G$, 
\[
C_{\sigma(G)}(\nu(G)) = \zeta_1(\nu(G))
\]

\noindent{\rm (e)}
Let $\Delta$ be one of the relations $\sn, \mathrm{asc}$, and let $K \lhd G$
with $K \in \acute{\EEE}\fA$. If $H \lhd G$ and 
$H$ contains every $\Delta\fN_2$ normal subgroup of $G$, then
\[
K \cap C_G(H) \subseteq C_K(H\cap K) =  \zeta_1(H\cap K)
\]
\end{lemma}
\begin{proof}
(a) This is in Robinson \cite[Lemma 2.19.1]{R68}. The proof is short so we include it.

Let $G \in \fZ$ with upper central series $\{Z_\alpha: 0 \leq \alpha \leq \tau\}$.
Let $A$ be a maximal abelian normal subgroup of $G$ and let $C =C_G(A)$. Then
$C \lhd G$ and $A \subseteq C$. If $C \neq A$, we can find an ordinal
$\alpha$ minimal with respect to $C \cap Z_\alpha \not\subseteq A$.
Clearly $\alpha$ is not a limit ordinal. Let $x \in C \cap Z_\alpha \setminus A$.
If $B = \langle x \rangle A$ then $B' =\ONE$ and $[B,G] \subseteq [C \cap Z_\alpha,A]
\subseteq (C \cap Z_{\alpha-1})A \subseteq B$, so $B \lhd G$, which contradicts maximality of $A$.
So $C = A$.

(b) Evidently if $M \in \acute{\EEE}\fA$ and $H \leq M$, then
$H \in \acute{\EEE}\fA$. Since $\sigma(G) \subseteq K$ it is enough to show that 
$K \in \acute{\EEE}\fA$. We may assume $K \neq \ONE$, and define $K_0 = \ONE$.

Suppose that $K_\alpha$ has been defined so that $K_\alpha \lhd G$ and
$K_\alpha \subseteq K$. Then there is $H \lhd G$ such that
$H \in \fZ \cup \EEE\fA$ and $H \not\subseteq K$. 
There are now two cases.

If $H \in \fZ$, let  $\beta$ be minimal with respect to $\zeta_\beta(H) = H_\beta \leq K_\alpha$.
Then $\beta$ is not a limit ordinal so $\beta-1$ exists and
$\zeta_{\beta-1}(H) \subseteq H_\beta$. Put $K_{\alpha+1} = K_\alpha H_\beta$, 
so that $K_{\alpha+1}' \subseteq K_\alpha H_\alpha \subseteq K_\alpha\zeta_{\beta-1}(H)
\subseteq K_\alpha$, since $K_\alpha \lhd G$. Furthermore, $\zeta_\beta(H)\, \mathrm{ch}\,
H \lhd G$, so $\zeta_\beta(H) \lhd G$ and $K_{\alpha+1} \lhd G$.

If, on the other hand, $K \in \EEE\fA$, then $H^{(d)} = \ONE$ for some $d$.
There exists $n \leq d$, maximal with respect to $H^{(n)} \not\subseteq K_\alpha$.
Then $H^{(n+1)} \subseteq K_\alpha$. In this case put
$K_{\alpha+1} = K_\alpha H^{(n)}$, so that $K_{\alpha+1}' \subseteq K_\alpha H^{(n+1)} 
\subseteq K_\alpha$. Further, $H^{(n)}\, \mathrm{ch}\, H \lhd G$, and $K_{\alpha+1} \lhd G$.

In either case, we see that if $K_\alpha \subseteq K$, $K_\alpha \lhd G$,
and $K_\alpha \neq K$, we can define $K_{\alpha+1}$ with
$K_\alpha < K_{\alpha+1}$, $K_{\alpha+1} \lhd G$, and $K_{\alpha+1}' \subseteq K_\alpha$.
If $\tau$ is a limit ordinal and $K_\alpha$ has been defined for all
ordinals $\alpha < \tau$, put $K_\tau = \bigcup_{\alpha < \tau} K_\alpha$.
By transfinite induction we can continue this process until some $K_\mu = K$.
Hence $K$ has an ascending series of normal subgroups of $G$ with abelian factors.
In particular $K \in \acute{\EEE}\fA$ so $\sigma(G) \in\acute{\EEE}\fA$.

(c) The result is obvious if $\Delta$ is $\leq$, for then $H = G$. Let $\Delta$
be one of $\lhd^\lambda$ ($\lambda$ an ordinal $> 0$), $\sn$, $\mathrm{asc}$,
and let $\{G_\alpha : 0 \leq \alpha \leq \tau\}$ be an ascending series of
$T$ with $G_\alpha \lhd G$ for all $\alpha \leq \tau$. Let $C = C_G(H)$.
Then $C \lhd G$, so $C \cap G_\alpha \lhd G_\alpha$ for all $\alpha \leq \tau$.
Suppose that $C \not \subseteq H$. Then we can find $\alpha$ minimal with
respect to $C\cap G_\alpha \not \subseteq H$. Evidently $\alpha > 0$ and
$\alpha$ is not a limit ordinal, so $\alpha-1$ exists. Now
\[
(C \cap G_\alpha)' \leq C \cap G_\alpha' \leq G \cap G_{\alpha-1} \leq H
\]
by minimality of $\alpha$, so 
\[
\gamma_3(C \cap G_\alpha) \leq [H,C] = 1
\]
whence $C\cap  G_\alpha \in \fN_2$, so $C\cap  G_\alpha \leq H$, a contradiction
since $C \cap G_\alpha \lhd G$. Thus $C \leq H$.

(d) Clearly $C = K \cap C_G(H) \lhd G$. First we show that
$C \leq H$, hence that $C \subseteq \zeta_1(H)$. By the proof of part (b) above,
$K$ has an ascending abelian series $\{K_\alpha :0 \leq \alpha \leq \mu\}$
such that $K_\alpha \lhd G$ for all $\alpha \leq \mu$. Suppose that $C \not \subseteq H$. 
Then we can find $\alpha$ minimal with
respect to $C\cap K_\alpha \not \subseteq H$. Evidently $\alpha > 0$ and
$\alpha$ is not a limit ordinal, so $\alpha-1$ exists. Since $K_\alpha' \subseteq K_{\alpha-1}$,
we have $(C\cap K_\alpha)' \subseteq C\cap K_{\alpha-1} \subseteq H$,
whence $C\cap  K_\alpha \in \fN_2$ and $C\cap  K_\alpha \lhd G$, so
$C\cap  K_\alpha \leq H$, a contradiction. Therefore $C \leq H$, so $C \leq \zeta_1(H)$.
But $\zeta_1(H)$ is an abelian normal subgroup of $G$, so 
$\zeta_1(H) \leq K\cap C_G(H) = C$. Thus $C = \zeta_1(H)$.
Evidently
\[
C \leq C_K(H \cap K) \cap H = C_{H\cap K}(H\cap K) = \zeta_1(H\cap K)
\]
In particular, $R = \sigma(G) \leq K$, so if we put $H = \nu(G)$,
then
\[
C_R(H) = R \cap C_G(H) \leq \zeta_1(H\cap K) = \zeta_1(H)
\]
since $H \leq R \leq K$. Conversely,
\[
\zeta_1(H) \leq R \cap C_G(H) = C_R(H)
\]

(e) Trivially $C_G(H) \leq C_G(H \cap K)$, so $K \cap C_G(H) \leq C_K(H \cap K)$.
Let $N \in \fN_2$ and $N \lhd K$. Since $K \Delta G$, we have $N \Delta G$,
so $N \leq H$, whence $N \leq H \cap K$ and $H \cap K$ contains every
$\Delta\fN_2$-subgroup of $K$. Further, $H \cap K \lhd K$, so by part (c)
$C_K(H\cap K) \leq H\cap K$.
\end{proof}

\begin{corollary}
\label{C:MC1.3}
Let $G$ be a group. Then

\noindent{\rm (a)}
$\nu(G) \in \fF \cap \fN$ if and only if $\sigma(G) \in \fF \cap \EEE \fA$.

\noindent{\rm (b)}
If $G \in {\rm Max}\mbox{-}\lhd\fA$ and $A$ is an abelian normal subgroup of $G$,
then to each $\fZ$-normal subgroup  $K$ of $G$ there corresponds an integer $n$
such that $A \cap K = A \cap \zeta_n(K)$.


\end{corollary}

\begin{proof}

(a) Let $H = \nu(G) \in \fF \cap \fN$ and let $R = \sigma(G)$.
Then $R/C_R(H) \in \fF$, and by Lemma~\ref{L:MC1.2}(d) we have $C_R(H)
= \zeta_1(H) \in \fF$ so $R \in \fF$. Clearly $R \in \EEE\fA$. The converse is trivial
since $\nu(G) \in \LLL\fN$ and $\nu(G) \subseteq \sigma(G)$.

(b) For each $\alpha$, we have $\zeta_\alpha(K)\,\mathrm{ch}\,K \lhd G$, so
$\zeta_\alpha(K) \lhd G$, whence $A \cap \zeta_\alpha(K)$ is 
an abelian normal subgroup of $G$. Since $G \in {\rm Max}\mbox{-}\lhd\fA$ we
can find $n$ such that $A \cap \zeta_n(K) = A \cap \zeta_{n+1}(K) = \cdots$,
so $A \cap \zeta_\omega(K) = A \cap \zeta_n(K)$. If $A \cap \zeta_\alpha(K)
\leq A \cap \zeta_n(K$, then
\[
[ A \cap \zeta_{n+1}(K),K] \subseteq  A \cap \zeta_\omega(K) \subseteq  A \cap \zeta_n(K)
\]
so $A \cap \zeta_{\alpha+1}(K)\subseteq A \cap \zeta_{n+1}(K) = A \cap \zeta_n(K)$
If $\lambda$ is a limit ordinal and $A \cap \zeta_\alpha(K) = A \cap \zeta_n(K)$
for all $\alpha < \lambda$, then $A \cap \zeta_\lambda(K) 
= \bigcup_{\alpha < \lambda} (A \cap \zeta_\alpha(K)) \subseteq A \cap \zeta_n(K)$.
Therefore $A \cap K = A \cap \zeta_n(K)$.


\end{proof}

\ignore{
We do not know whether,
if $G \in \fZ \cap {\rm Max}\mbox{-}\lhd\fA$, then $G \in \fG \cap \fN$.

Indeed, we do not even know whether,
if $G \in \fN \cap {\rm Max}\mbox{-}\lhd\fA$, then $G \in \fG$.
}

\section{Hypercentral Groups in Max-$\lhd \fA$ are\\Finitely Generated Nilpotent}
\label{S:HGMFGN}

In this section we prove that every hypercentral group 
that lies in Max-$\lhd \fA$ is finitely generated nilpotent. 

Recall that the Hirsch--Plotkin radical $\rho(G)$ is the unique maximal
locally nilpotent normal subgroup of $G$. The {\em upper Hirsch--Plotkin series}
$\rho_\alpha(G)$ of $G$ for ordinals $\alpha$ is defined recursively by
\beqn
\rho_1(G) &=& \rho(G) \\
\rho_{\alpha+1}(G)/\rho_\alpha(G) &=& \rho(G/\rho_\alpha(G)) \\
\rho_\lambda(G) &=& \cup_{\alpha<\lambda} \rho_\alpha(G)
	\qquad \lambda\ \mbox{a limit ordinal}
\eeqn

The next definition is from Robinson \cite{R70a} page 59.
\begin{definition}\em
\label{D:radical_group}
A {\em radical group} is a group that is {\em hyper-locally-nilpotent}.
That is, the upper Hirsch-Plotkin series reaches $G$
after some ordinal number of steps: $\rho_\alpha(G)=G$.

In particular, every hypercentral group is radical. 
\end{definition}

Next, we cite a result of Robinson \cite{R68}, generalising
theorems of Mal'cev \cite[Theorem 2]{M51}, Smirnov \cite[Theorem 3]{S53}, and Baer \cite[Satz B$'$]{B55}:

\begin{lemma}
\label{L:radical_group}
Any radical subgroup of the automorphism group of a polycyclic-by-finite 
group is polycyclic.
\end{lemma}
\begin{proof}
See Robinson \cite[Theorem 3.27]{R68}.
\end{proof}

Define $\aut(G)$ to be the automorphism group of $G$.

\begin{theorem}
\label{T:ELNmax-An}
Let $G$ be a radical group. 
Suppose that $G$
has a self-centralising abelian normal subgroup that is finitely generated.
Then $G$ is finitely generated nilpotent.
\end{theorem}
\begin{proof}
Let $A$ be a finitely generated abelian normal subgroup
with $C_G(A) = A$.
Then $H =G/A$ acts faithfully on $A$  by conjugation, so $H$ embeds
in $\aut(A)$.
Since $A \in \fA \cap \fG$ it is in particular polycyclic.                                     
Lemma \ref{L:radical_group} implies that $H$ is polycyclic.
In particular, $H \in \fG$. 
Now $A \in \fG$ and $G/A = H \in \fG$, so by $\EEE$-closure
$G \in \fG$. Since it is hypercentral, hence locally nilpotent,
$G \in \fN\cap\fG$.
\end{proof}

We easily deduce:

\begin{theorem}
\label{T:ZAmax-An}
Let $G$ be hypercentral and in ${\rm Max}-\lhd\fA$. Then $G$ is finitely generated nilpotent.
\end{theorem}

\begin{proof}
Let $A$ be a maximal abelian normal subgroup of $G$. By Lemma \ref{C:MC1.3}(b),
with $K = G$, there is an integer $n$ such that $A \cap \zeta_n(G) = A$.

Define $A_i = A \cap \zeta_i(G)$. Then there is a chain
\[
A_0 \subseteq A_1 \subseteq A_2 \subseteq \cdots \subseteq A_n = A
\]
of abelian normal subgroups of $G$.

Suppose that $B$ is a subgroup such that $A_i \subseteq B \subseteq A_{i+1}$.
Then $[B,G] \subseteq [A,G]  \subseteq A$. Also $[B,G] \subseteq [\zeta_{i+1}(G),G]  \subseteq \zeta_i(G)$. Therefore $[B,G] \subseteq A \cap \zeta_i(G) = A_i \subseteq B$. Thus
$B$ is a normal abelian subgroup of $G$.

It follows that $A_{i+1}/A_i \in \fG$ since $A_{i+1}/A_i \in \fA \cap {\rm Max}$.
Therefore by $\EEE$-closure of $\fG$ we have $A_i \in \fG$ for all $i$.
But $A = A_n$ so $A \in \fG$.

We know that  $C_G(A)=A$ by Lemma \ref{L:MC1.2}(a).
The group $G$ is hypercentral, hence locally nilpotent, hence radical.
Now Theorem \ref{T:ELNmax-An} implies that
$G \in \fN\cap\fG$.
\end{proof}

\begin{remark}\em
(a) Therefore $G$ satisfies ${\rm Max}$, because $\fN\cap\fG \subseteq {\rm Max}$
 by Baer \cite{B55}.

(b) When $G \in \fN$, or indeed when $G \in \EEE\fA$,  the
weaker theorem of Mal'cev \cite{M51} (replacing 'radical' by `polycyclic') suffices in the proof,
in place of Lemma \ref{L:radical_group}.

(c) Conversely, if $G \in \fN\cap\fG$ then $G \in {\rm Max}$.
This easily implies the converses of Theorems \ref{T:ELNmax-An} and \ref{T:ZAmax-An}.
\end{remark}


\section{Maximal Conditions for Abelian Normal Subgroups of Quotients}

We relate the class $\mda$ to $\qn$ and ${\rm Max}-\lhd$.

For any class of groups $\fX$ define
\[
\fX^{\QQQ} = \{G : G/N \in \fX\ \forall N\lhd G \}
\]
\begin{lemma}
\label{L:QN2.1}
For all $k \geq 1$, 
\[
(\qn\fA)^{\QQQ} = (\qn\fA^k)^{\QQQ}
\]
\end{lemma}

\begin{proof}
Clearly ${\rm Max}-\lhd\fA^k \subseteq {\rm Max}-\lhd\fA$,
so $({\rm Max}-\lhd\fA^k)^{\QQQ} \subseteq ({\rm Max}-\lhd\fA)^{\QQQ}$.

To prove the reverse inclusion
\begin{equation}
\label{e:maxAk}
({\rm Max}-\lhd\fA)^{\QQQ} \subseteq ({\rm Max}-\lhd\fA^k)^{\QQQ}
\end{equation}
we use induction on $k$. Clearly~\eqref{e:maxAk} is valid for $k = 1$.
So we assume it holds for $k$ and prove it for $k+1$.

Let
\begin{equation}
\label{e:acc2}
N_0 \subseteq N_1 \subseteq N_2  \subseteq \cdots 
\end{equation}
be an increasing chain of $\fA^{k+1}$ normal subgroups of $G$. For each $m \in \N$,
\[
N_m^{(k)} \ \mbox{ch}\ N_m \lhd G
\]
so $N_m^{(k)} \lhd L$. Further, $N_m^{(k)} \in \fA$. So
\[
N_0^{(k)} \subseteq N_1^{(k)} \subseteq N_2^{(k)}  \subseteq \cdots 
\]
is an ascending chain of abelian normal subgroups of $G \in ({\rm Max}-\lhd\fA)^{\QQQ}$.
By ${\rm Max}-\lhd\fA$ it terminates, so there exists $r \in \N$ such that
\[
N_m^{(k)} = N_r^{(k)} \quad \mbox{for all}\ m \geq r
\]
Now
\[
N_r/N_r^{(k)} \subseteq N_{r+1}/N_r^{(k)} \subseteq N_{r+2}/N_r^{(k)}  \subseteq \cdots 
\]
is an ascending chain of $\fA^k$ normal subgroups of $G/N_r^{(k)}$. Since by induction
$G \in ({\rm Max}-\lhd\fA)^{\QQQ} = ({\rm Max}-\lhd\fA^k)^{\QQQ}$, the chain
terminates. Therefore~\eqref{e:acc2} terminates.
\end{proof}

We then have a characterisation of Noetherian groups in terms
of quasi-Noetherian groups and the property $({\rm Max}-\lhd\fA)^{\QQQ}$:

\begin{theorem}
\label{T:qnmaq}
\[
{\rm Max}-\lhd = ({\rm Max}-\lhd\fA)^{\QQQ} \cap \qn 
\]
\end{theorem}

\begin{proof}
Clearly
\[
{\rm Max}-\lhd \subseteq ({\rm Max}-\lhd\fA)^{\QQQ} \cap \qn 
\]
To prove the reverse inclusion, let $G \in ({\rm Max}-\lhd\fA)^{\QQQ} \cap \qn$,
and consider an increasing chain of normal subgroups of $G$:
\begin{equation}
\label{e:acc3}
N_0 \subseteq N_1 \subseteq N_2  \subseteq \cdots \quad \mbox{where}\ \bigcup_n N_n = N
\end{equation}
For a contradiction, assume that $N_n \neq N_{n+1}$ for
all $n \in \N$. Since $G \in \qn$ there exists $m \in \N$ such that
\[
[G^{(m)},N] \subseteq N_m
\]
Therefore $N^{(m+1)} \subseteq I_m$, so $N/N_m \in \fA^{m+1}$.
So
\begin{equation}
\label{e:acc4}
\frac{N_m}{N_m} \subseteq \frac{N_{m+1}}{N_m} \subseteq \frac{N_{m+2}}{N_m} \subseteq \cdots
\end{equation}
is an ascending chain of normal subgroups of $G/N_m$.

Since $G \in ({\rm Max}-\lhd\fA)^{\QQQ}$, we have $G/N_m \in ({\rm Max}-\lhd\fA^{k+1})^{\QQQ}$
by Lemma~\ref{L:QN2.1}. Therefore~\eqref{e:acc4} terminates. So~\eqref{e:acc3} terminates,
and $G \in {\rm Max}-\lhd$ as claimed.
\end{proof}




\end{document}